\documentclass{article}

\usepackage[utf8]{inputenc}
\usepackage[english]{babel}
\usepackage[utf8]{inputenc}
\usepackage{amssymb}
\usepackage{amsmath}
\usepackage{amsthm}
\usepackage{fancyhdr}
\usepackage{changepage}
\usepackage{mathtools}
\usepackage{multicol}
\usepackage{wrapfig}
\usepackage{xcolor}
\usepackage{enumitem}
\usepackage{hyperref}
\usepackage{mathrsfs}
\usepackage{natbib}
\usepackage{thmtools}
\usepackage{nameref}

\providecommand{\dotdiv}{
  \mathbin{
    \vphantom{+}
    \text{
      \mathsurround=0pt
      \ooalign{
        \noalign{\kern-.35ex}
        \hidewidth$\smash{\cdot}$\hidewidth\cr 
        \noalign{\kern.35ex}
        $-$\cr
      }
    }
  }
}
\usepackage[margin=1.5in]{geometry}
\newtheorem{thm}{Theorem}[section]

\newtheorem{lem}[thm]{Lemma}
\theoremstyle{definition}
\newtheorem{dfn}[thm]{Definition}
\theoremstyle{remark}

\title{Arrow's Impossibility Theorem: \linebreak Computability in Social Choice Theory}
\author{Alex Hall TEWH}
\date{May 2022}

\begin{document}

\maketitle

\section{Intuition Behind Arrow's Theorem}

Arrow's impossibility theorem is a classic result in social choice theory. It establishes that a certain group of desirable conditions for aggregating voters' preferences into a single ``verdict" cannot all simultaneously hold. These are summarised by \nameref{dfn:arrow}. A voter profile is a weak ordering which essentially represents a voter's ranked preference with respect to a list of alternatives, allowing for the possibility of ties\footnote{It is worth remarking that we can disallow ties if preferred, as this extension actually makes Arrow's result stronger.}. A social welfare function is then a map which takes the entire voting population's preferences and outputs a single profile. In an applied setting we can think of this as a verdict which a governing body can act upon, for example by appealing to some majority rule which determines the orderings. Whilst any valid function is a social welfare function, it is typically desirable to have a degree of ``fairness" in the outcome. Arrow formalised three fairness conditions within the rest of his axioms. First, representability is desired for the welfare function, so at a minimum if voters unanimously agree that alternative A is better than alternative B, then the social welfare function must respect this and prefer A over B. Second, the outcome should consider alternatives independently. This means that if we look at two distinct votes in which the voters feel identically about candidates A and B, then even if they differ on their preferences with respect to candidate C, the relative placement of A and B must not change. Third, in order for the vote to have any meaning, it should not be dictated by a single voter. This condition captures the very reason for voting; after all, a vote which is dictated by a single voter is not worth participating in for all other voters. It is worth emphasising that the dictator in this scenario is the same regardless of other voters' inputs for a fixed social welfare function. The significance of Arrow's theorem is really in the fact that the first two conditions exclude the third, meaning that the two seemingly reasonable demands for unanimity and independence undermine democracy by ensuring dictatorship. \par
It is helpful to consider a small example taken from \cite{sep-arrows-theorem}
on the SEP, to see how we can encounter difficulties in creating social welfare functions. Consider this scenario with three voters and candidates A, B and C. We use $X > Y$ to express that $X$ is preferred over $Y$:
\begin{enumerate}
    \item A $>$ B $>$ C
    \item C $>$ A $>$ B
    \item B $>$ C $>$ A
\end{enumerate}
The example is generated by shifting the preferences for the three voters, so that no agreement seems possible. A seemingly intuitive pairwise majority rule doesn't work as it produces a cyclic result that ``A $>$ B $>$ C $>$ A $>$ ..." which clearly isn't possible. The next ``fairest" thing to do might be to consider the result a three-way tie. Since our ordering is weak, this is permissible. The problem is that we've only reached a verdict by considering the preferences on an \textit{ad hoc} basis. A social welfare function on the other hand must be fixed before the voter's preferences are input, capturing the desire that the outcome is determined impartiality. Any attempt to generalise rules which lead to the outcome ``A = B = C" fail because of Arrow's theorem. For example, if we decide all candidates draw unless the vote is unanimous on all inputs, we fail to meet independence. On the other hard if we maintain independence by always ruling that all candidates draw and neglect the votes, we fail to respect a unanimous ruling. Finally, always listening to some fixed candidate satisfies the first two demands but is just a dictatorship. \par

\section{Goals of the Paper}

Since Arrow first published his theorem, there have been significant developments. \citeauthor{fishburn1970arrow} extended Arrow's analysis to infinite sets of voters and was able to show that in this case, it was actually possible to simultaneously satisfy all the axioms, meaning non-dictatorial social welfare functions did exist. \citeauthor{mihara1995arrow} suggested that this infinite extension is not just a toy result of mathematics but has genuine applications over finite sets of voters if the voters are allowed to express countably many different preferences contingent on some set of future states. The vote is modelled by treating each original voter under each state as an independent voter which results in an infinite case. In the meantime, \citeauthor{kirman1972arrow} re-proved Arrow's and Fishburn's results using the existing mathematical concept of ``ultrafilters", taken from topology and set theory. Later, \citeauthor{lewis1988infinite}, \citeauthor{mihara1995arrow} and \citeauthor{benedict2022arrow} interpreted the results of Arrow's theorem through the lens of computability. \par
In this paper we will be surveying the results of \cite{mihara1995arrow, mihara1997arrow, mihara1999arrow}. We aim to differ from the surrounding literature by focusing on mathematical and computability theoretic content rather than applied results. In doing this, we are able to be more comprehensive, and take time to make the results approachable for those with only a general background in mathematics or logic. This approach should combine the best of both \cite{kirman1972arrow} and \cite{mihara1997arrow} in a more concise form, and with additional explanation to help build an intuitive understanding of the results. In particular, borrowing from the approach of \cite{benedict2022arrow} we place more importance on the logical groundings of computability theory. \citeauthor{mihara1995arrow} himself spends very little time on this even though it underpins the motivation for his main result. We also end with a section that introduces advanced concepts of relative computability. extending beyond what is currently found in the referenced papers. We suggest ways in which these advanced concepts could lead to further research.

To achieve all this, we begin by introducing the classical formulation of \nameref{dfn:arrow} by \cite{kirman1972arrow}. We will introduce sufficient notation so that readers without any background knowledge can interpret the statements of \nameref{thm:arrowsfinite} and \nameref{thm:fishburn}. We will show some proofs leading up to but excluding the main result, because they go into significant technical depth. By explaining the intuition regarding ultrafilters, as well as proving a few of their main properties such as \ref{lem:ultrafinitefixed}, \ref{lem:ultracomplement} and \ref{lem:singletongenerated}, the original proofs should be made significantly more accessible. At this point we introduce a significant amount of background on computability theory, again in sufficient detail that readers should require no specific background knowledge. We then have a cursory introduction on second-order arithmetic. Combining this we get computability, which roughly describes the sorts of things we can compute by deterministic means, and arithmetic, which defines a method of expression within rigorous bounds. This means that our computations can then be expressed in the language of this arithmetic. \par
This is where \citeauthor{mihara1995arrow} contributes, and we primarily follow. Because in \cite{fishburn1970arrow}, \nameref{thm:fishburn} uses a non-constructive method\footnote{i.e. without giving an explicit example.} to show the existence of his non-dictatorial social welfare functions, it is natural to wonder what this sorts of functions would look like. Short of finding a constructive proof, the most reasonable step is to attempt to restrict the properties of these non-dictatorial functions. Unfortunately, this infinite case does not lend itself to coding for use in computability; instead \citeauthor{mihara1995arrow} proposed his own modified concept of pairwise computability that demands we have some computable method for determining the relative positions of any pair of candidates in the verdict. We will leave discussion of the consequence of pairwise computability until later when the methods are better understood. His final result, as is ours, is then to prove that even with infinitely many voters like in \citeauthor{fishburn1970arrow}'s case, we must have all pairwise computable social welfare functions being dictatorial, which undermines any hope we might have gained by rejecting Arrow's theorem.

\section{Building to Arrow's Classical Result}

Before we can define \nameref{dfn:arrow}, we first formally define the notation and specifics of the concepts previously introduced. To do this, we need to understand orderings and the rules that govern them. We need to introduce the notation used for voters and the alternatives they vote on. Finally we need to understand the way in which our social welfare function is seen as a ``function".

\begin{dfn}
An order $P$ on a set $X$ is a binary relation on the elements of $X$. We express the order relation $P$ holding on elements $x,y \in X$ by writing $xPy$. An order is called a weak order if $\forall x,y,z \in X$ it is both
\begin{enumerate}[label=(O\arabic*)]
    \item\label{dfn:weakorderasym} asymmetric, meaning $xPy \implies \neg yPx$ and 
    \item\label{dfn:weakordernegtrans} negatively transitive, meaning $( \neg xPy \land \neg yPz ) \implies \neg xPz$.
\end{enumerate}
\end{dfn}

\begin{lem}
For a weak order $P$ on $X$ and $\forall x,y,z \in X$ we have that
\begin{enumerate}[label=(\roman*)]
    \item\label{lem:weakorder1} $P$ is transitively incomparable, meaning $xPz \implies ( xPy \lor yPz )$ and
    \item\label{lem:weakorder2} $P$ is transitive, meaning that $( xPy \land yPz ) \implies xPz$ and
    \item\label{lem:weakorder3} $( ( \neg yPx \land yPz ) \lor ( xPy \land \neg zPy ) ) \implies xPz$.
\end{enumerate}
\end{lem}
\begin{proof}
We begin by proving \ref{lem:weakorder1} using the contrapositive of \ref{dfn:weakordernegtrans} and De Morgan's Law which gives 
$$xPz \implies \neg ( \neg xPy \land \neg yPz ) \implies ( xPy \lor yPz ).$$
We can see that \ref{lem:weakorder2} arises now by using \ref{lem:weakorder1} since 
$$( xPy \land yPz ) \implies ( ( xPz \lor zPy ) \land yPz ).$$ 
Using property \ref{dfn:weakorderasym} we now further get 
$$( xPy \land yPz ) \implies ( ( xPz \lor zPy ) \land \neg zPy ) \implies xPz$$
as desired, since assuming $zPy$ in the ``OR" clause gives an immediate contradiction. \\
To prove \ref{lem:weakorder3} we assume $\neg xPz$ for a contradiction. We proceed by cases analysis, so if we assume the first clause we get
$$ ( \neg yPx \land \neg xPz \land yPz ) \implies ( \neg yPz \land yPz ) \implies \bot $$
and if we assume the second clause we get
$$ ( xPy \land \neg xPz \land \neg zPy ) \implies ( xPy \land \neg xPy ) \implies \bot $$
where in both cases we've simply used \ref{dfn:weakordernegtrans}. Thus we have $xPz$ as desired.
\end{proof}

\begin{dfn}
If $V$ is a set of voters, $X$ is a set of alternatives and $W = \{ P \subseteq X \times X : P \text{ is a weak order on } X \}$ is the set of all weak orders on $X$, then we call any function $f: V \to W$ a ``voter profile" which maps each voter to their preferred order, meaning $xf(v)y$ represents voter $v$ preferring $y$ to $x$. We further let $F = \{ f : V \to W\}$ be the set of all profiles.
\end{dfn}

\begin{dfn}
A function $\sigma$ is called a ``social welfare function" if it maps any profile to a single weak order, i.e. if $\sigma : F \to W$. The order corresponding to a fixed profile can be thought of as the ``verdict" our function decides w.r.t. the voters' choices. 
\end{dfn}

Intuitively, each profile describes a possible way that the voters can rank alternatives. In an applied setting, the outcome of an election can be seen as the single realised profile. Clearly however, before an election is run we cannot know which profile will be realised without rigging the outcome. As such, if we want to agree on a method to determine the verdict beforehand, we must be able to map every possible profile to a weak order which is the outcome.

\begin{dfn}
If $a,b \in X$ and $U \subseteq V$ and $f$ is a profile on $X$ and $V$, then we write $af(U)b$ to denote that all voters on $U$ prefer $b$ to $a$, i.e. that $\forall u \in U, af(u)b$. \\
We further write that profiles $f(u) = g(u)$ on $\{a,b\}$ whenever their preferences for a given voter $u$ agree on $a$ and $b$ specifically, i.e. when $( af(u)b \iff ag(u)b ) \lor ( bf(u)a \iff bg(u)a )$. \\
Finally we write $f = g$ on $\{a,b\}$ whenever their preferences agree on $a$ and $b$ for all voters, i.e. $\forall v \in V, f(v) = g(v)$ on $\{ a,b \}$.
\end{dfn}

At this point, we have defined sufficient notation to fully express the statements of Arrow's Axioms formally.

\begin{dfn}[Arrow's Axioms]\label{dfn:arrow} \quad
\begin{enumerate}[label=(A\arabic*)]
    \item\label{A1:alternatives}[Alternatives] There are at least $3$ alternatives; so $\lvert X \lvert \geq 3$.
    \item\label{A2:welfarefunction}[Welfare Function] $\sigma$ is a social welfare function on alternatives $X$ and voters $V$.
    \item\label{A3:unanimity}[Unanimity] If ever all voters agree on an order this must be respected; so 
    $$\forall a,b\in X, \forall f \in F, af(V)b \implies a\sigma(f)b.$$
    \item\label{A4:independence}[Independence] Only preferences on the ordering of $a$ and $b$ should affect
    our verdict of $a$ relative to $b$, meaning functions which are identical on $\{ a,b \}$ have the same verdict; so
    $$\forall a,b \in X, \forall f,g \in F, ( f = g \text{ on } \{a,b\} ) \implies  ( \sigma(f) = \sigma(g) \text{ on } \{a,b\} ).$$
    \item\label{A5:nondictatorial}[Non-dictatorial] No one voter should dictate the verdict in all profiles; so
    $$\neg \exists v_0 \in V \text{ s.t. } \forall a,b \in X, \forall f \in F, af(v_0)b \implies a\sigma(f)b.$$
\end{enumerate}
When axiom \ref{A5:nondictatorial} is violated, we say that the system is dictatorial.
\end{dfn}

As alluded to, the classical result of Arrow's Theorem is to show that all social welfare functions over finite voter sets are dictatorial. We follow the proof of \cite{kirman1972arrow} to achieve this. They take an approach which involves likening social welfare functions to topological/set-theoretic objects called (ultra)filters. This allows us to work abstractly and re-use some of the established properties of these objects, rather attempting to reason about the social welfare functions directly.

\begin{dfn}
Given a set $V$, a system of subsets $\mathcal{F}$ is called a filter on $V$ if $\forall A,B \subseteq V$ it is
\begin{enumerate}[label=(F\arabic*)]
    \item\label{dfn:upwardclosure} upward closed, meaning $( A \subseteq B \land A \in \mathcal{F} ) \implies B \in \mathcal{F}$ and
    \item\label{dfn:finiteintersection} closed under (finite) intersection, meaning $A,B \in \mathcal{F} \implies A \cap B \in \mathcal{F} $ and
    \item\label{dfn:proper} proper, meaning $\emptyset \notin \mathcal{F}$.
\end{enumerate}
\end{dfn}

\begin{dfn}
A filter $\mathcal{F}$ is called ``free" (or non-principal) if $\bigcap\mathcal{F} = \emptyset$. If a filter isn't free it is called ``fixed" (or principal).
\end{dfn}

\begin{lem}\label{lem:ultrafinitefixed}
If $V$ is a finite set and $\mathcal{F}$ is a filter on $V$, then $\mathcal{F}$ is fixed.
\end{lem}
\begin{proof}
Since $V$ is finite, we can use finite induction on \ref{dfn:finiteintersection} to show that $\bigcap\mathcal{F} \in \mathcal{F}$, where we've used the fact that $V$ being finite also implies it has finitely many subsets, and thus that the algebra and filter are also finite. Then clearly if $\bigcap\mathcal{F} = \emptyset \in \mathcal{F}$, we violate \ref{dfn:proper} for a contradiction. So $\mathcal{F}$ is not free and is fixed.
\end{proof}

\begin{dfn}
If $\mathcal{F}$ and $\mathcal{G}$ are filters on $V$, we say $\mathcal{F}$ is strictly finer than $\mathcal{G}$ when $\mathcal{G} \subsetneq \mathcal{F}$.
\end{dfn}

\begin{dfn}
If $\mathcal{U}$ is a filter on $V$, we say it is an ultrafilter whenever $\neg \exists \mathcal{F}$ s.t. $\mathcal{F}$ is a strictly finer filter than $\mathcal{U}$ on $V$.
\end{dfn}

The definition of what an ultrafilter is does differ from text to text, with some defining it as a filter with the property that every set or its complement is included. However following \cite{kirman1972arrow}, we use the above definition. Instead the equivalent characterisation can be proven as a lemma, although this is not shown in \citeauthor{kirman1972arrow}. They instead introduce the concept of a filterbase and shown a similar lemma is true for filterbases at \cite[p.275]{kirman1972arrow}. To avoid having to define this concept, we show this result directly by an original proof.

\begin{lem}\label{lem:ultracomplement}
$\mathcal{U}$ is an ultrafilter on $V$ if and only if $\mathcal{U}$ is a filter and $\forall A \subseteq V$ either $A \in \mathcal{U}$ (exclusively) or $A^c = V \setminus A  \in \mathcal{U}$.
\end{lem}
\begin{proof}
Suppose $\mathcal{U}$ is a filter and $\forall A \subseteq V$, either $A \in \mathcal{U}$ or $A^c \in \mathcal{U}$. Suppose $\mathcal{U}$ is not maximal, and that $\mathcal{F}$ is a filter strictly finer than $\mathcal{U}$. Pick $X \in \mathcal{F} \setminus \mathcal{U} \neq \emptyset$. Since $X \subseteq V$ and $X \notin \mathcal{U}$, this implies $X^c \in \mathcal{U}$ by our assumed property. But then since $\mathcal{F}$ is also assumed to be strictly finer, $X^c \in \mathcal{F}$, and so by \ref{dfn:finiteintersection} we have $X \cap X^c = \emptyset \in \mathcal{F}$, but this is a contradiction with \ref{dfn:proper}. So no such $\mathcal{F}$ exists and $\mathcal{U}$ is an ultrafilter. \par
Now suppose $\mathcal{U}$ is an ultrafilter. Suppose that $B \subseteq V$ and $B \notin \mathcal{U}$. Assume $B \neq \emptyset$ as we already know $\emptyset^c = V \in \mathcal{U}$ by \ref{dfn:upwardclosure} for a non-trivial filter. Now let $\mathcal{F}$ be the upward closure (property \ref{dfn:upwardclosure}) of $\mathcal{U} \cup B$. Since $\mathcal{U}$ is an ultrafilter, $\mathcal{F}$ can't be a filter. But by construction \ref{dfn:upwardclosure} and \ref{dfn:proper} are satisfied, so we must violate \ref{dfn:finiteintersection} with some set in $\mathcal{F} \setminus \mathcal{U}$. WLOG we can see that we can choose $B$ to be this set as it is minimal in $\mathcal{F} \setminus \mathcal{U}$ by construction. So then $\exists C \in \mathcal{U}$ s.t. $B \cap C = D$ and $D \notin \mathcal{F}$. Suppose $D \neq \emptyset$, then just consider the set $\mathcal{F} \cup D$, which again can't be a filter (since it is a superset of $\mathcal{U}$) and repeat the argument. Instead therefore assume $B \cap C = \emptyset$. Then $C \subseteq B^c$, and since $C \in \mathcal{U}$, by \ref{dfn:upwardclosure} we get $B^c \in \mathcal{U}$. Therefore we've shown $\forall A \subseteq V$ either $A \in \mathcal{U}$ or $(A \notin \mathcal{U} \implies A^c \in \mathcal{U})$. To see that this is exclusive we only need to remark that $A^c, A \in \mathcal{U} \implies A \cap A^c = \emptyset \in \mathcal{U}$ by \ref{dfn:finiteintersection} which violates \ref{dfn:proper}.
\end{proof}

\begin{lem}\label{lem:singletongenerated}
If $V$ is a finite set, $\mathcal{U}$ is an ultrafilter on $V$ if and only if it has the form $\mathcal{C} = \{ U \subseteq V : v_0 \in U\}$ for a fixed $v_0 \in V$.
\end{lem}
\begin{proof}
First, in the reverse direction it is not hard to check that any set of the form $\mathcal{C}$ satisfies \ref{dfn:upwardclosure}, \ref{dfn:finiteintersection} and \ref{dfn:proper}, and thus is a filter. Clearly it is upward closed as any superset of a set containing $v_0$ still contains $v_0$ and thus is in $\mathcal{C}$. Since all finite intersections contain $v_0$ still, they are in $\mathcal{C}$. Finally $v_0 \notin \emptyset$. Now, since $\forall U \subseteq V$, either $v_0 \in U$ or (exclusively) $v_0 \in U^c$, we satisfy the conditions of \ref{lem:ultracomplement}, and so $\mathcal{C}$ is an ultrafilter. \par
In the other direction we suppose we have some unknown ultrafilter $\mathcal{U}$ over a finite set. By \ref{lem:ultrafinitefixed}, we know that $\bigcap\mathcal{U} \neq \emptyset$. So suppose $\bigcap\mathcal{U} = \{ v_0, v_1, \dots, v_n \}$ for distinct elements. First we suppose $n > 0$ for a contradiction. Clearly the singleton $\{ v_0 \}$ cannot be in the ultrafilter, else $v_1 \notin \bigcap\mathcal{U}$. But then using the forward direction of our proof, we know $\mathcal{C} = \{ U \subseteq V : v_0 \in U\}$ is a filter. But it is also strictly finer as it contains at least all of $\mathcal{U}$ as well as the singleton $\{ v_0 \}$. This contradicts the ultrafilter property, so we have that $n = 0$ exactly. In this case $\mathcal{U}$ takes exactly the desired form, as if the intersection is exactly $\{ v_0 \}$ we see all sets contain $v_0$.
\end{proof}

From the previous results, we can take away the intuition that an ultrafilter splits a set into an upward closed collection which are members of the filter, and their complements excluded from the filter. In the finite setting, we can also view them as ``generated"\footnote{This can be defined formally by repeatedly applying closure conditions required to have a filter.} by a single element $v_0$. From here we have sufficient machinery to understand the classical proof of Arrow's impossibility theorem given by \cite{kirman1972arrow}. We shall not see the specifics of the proof, as the \citeauthor{kirman1972arrow} paper is quite exhaustive on this matter. However we shall state the main theorems that lead to the conclusion without proof.

\begin{dfn}
Let $\Sigma$ be the set of all $\sigma$ which also satisfy \ref{A1:alternatives}, \ref{A2:welfarefunction}, \ref{A3:unanimity} and \ref{A4:independence}. Further let $\widetilde{V}$ be the set of all ultrafilters over $V$.
\end{dfn}

\begin{thm}\label{thm:ks1} Let $g : \Sigma \to \widetilde{V}$ be a mapping of $\sigma \mapsto \mathcal{U}_\sigma$.
\begin{enumerate}
    \item  We specify $\mathcal{U}_\sigma$ by demanding that $\forall U \in \mathcal{U}_\sigma, \forall x,y \in X, \forall f \in F, xf(U)y \implies x \sigma(f) y$. Further this mapping is injective, i.e. each $\sigma$ is mapped to a unique ultrafilter.
    \item The mapping $g$ is also a surjection of $\Sigma$ onto $\widetilde{V}$.
\end{enumerate}
\end{thm}

\begin{thm}\label{thm:ks2}
An element $\sigma \in \Sigma$ also satisfies \ref{A5:nondictatorial} if and only if the corresponding ultrafilter defined in \ref{thm:ks1}, $g(\sigma) = \mathcal{U}_\sigma$, is free.
\end{thm}

\begin{thm}[Arrow's Theorem]\label{thm:arrowsfinite}
If $V$ is finite, then a social welfare function $\sigma$ satisfying \ref{A1:alternatives} - \ref{A4:independence} implies $\neg$\ref{A5:nondictatorial}. In other words, the five axioms are inconsistent. 
\end{thm}
\begin{proof}
This final result drops out quickly. Since $V$ is finite, by \ref{lem:ultrafinitefixed} is it fixed and not free. Then finally by the negation of \ref{thm:ks2}, we get $\neg$\ref{A5:nondictatorial}.
\end{proof}

The intuitive takeaway from the Kirman and Sondermann proof in the finite case relates to our remark about ultrafilters being generated by a singleton $\{ v_0 \}$. If $\sigma$ is a social welfare function in $\Sigma$, then under the mapping $g$ in \ref{thm:ks1}, $\mathcal{U}_\sigma$ is generated by the dictatorial voter. We can see this voter is dictatorial, as the property in \ref{thm:ks1} demanded that the welfare function always agree with $v_0$, since it is in all $U \in \mathcal{U}_\sigma$. \par
Fortunately \ref{thm:ks2} makes no reference as to whether $V$ should have to be finite or not. Thanks to \cite[p. 270]{kirman1972arrow}, we see the existence of the axiom of choice is sufficient to prove the existence of free ultrafilters over infinite sets\footnote{For those familiar with choice, the intuition is that we use Zorn's lemma. The system of all cofinite sets on infinite set $V$, called $\mathscr{F}$, forms a filter, and is clearly with empty intersection. Fineness forms an ordering, with the $V$ acting as a natural upper bound on any chain. With respect to this ordering, being a maximally fine element is equivalent to being an ultrafilter. Thus the chain containing $\mathscr{F}$ contains a maximal element which is a free ultrafilter.}. This means that with the machinery we already have, we can prove the next result in the historical ``back and forth" of Arrow's theorem, which was Fishburn's Possibility Theorem, originally proven in \cite{fishburn1970arrow}.

\begin{thm}[Fishburn's Possibility Theorem]\label{thm:fishburn}
If $V$ is infinite, the axioms \ref{A1:alternatives} - \ref{A5:nondictatorial} are consistent.
\end{thm}
\begin{proof}
The result drops out of \ref{thm:ks2} exactly as in the proof of \nameref{thm:arrowsfinite} but with a free instead of fixed ultrafilter, giving the opposite conclusion.
\end{proof}

Classically, \nameref{thm:fishburn} is seen as giving some degree of ``hope" in the search for desirable voting systems, even if just in infinite cases. As previously discussed, \cite[p. 259-260]{mihara1997arrow} gives examples of practical uses for the infinite case in the real world setting. However in this paper we shall see that this hope is short-lived. We assume very gentle restrictions on what sorts of social welfare functions are permitted, viewed through the lens of computability theory. We then show that this restriction is already enough to lose ``possibility" in the infinite case, and we reaffirm that all social welfare functions (of our specification) are dictatorial. This shows the extreme limitations of Fishburn's positive result and implies using voting systems which can be described by our social welfare functions is still practically impossible without dictatorship.

\section{Introducing Computability Theory}

The main focus of \citeauthor{mihara1995arrow} across his three papers published on this topic was to consider whether we can actually use the free ultrafilters which are shown to exist. The Axiom of Choice, used in the proof of \nameref{thm:fishburn} is  non-constructive here, meaning that it affirms the existence of these desirable social welfare functions without providing us a method to obtain them. Computability Theory is a field of study which focuses on formalising the potential class of objects which we can obtain by ``algorithms", and is generally accepted as the limit of what is decidable unambiguously (or relative to a decision maker). Considering that we generally want our algorithms to be able to take inputs and return the same response every time, they are translated into the language of mathematics as functions\footnote{This is not they only way to see computability; we shall address Turing machines also.}. Thus in many ways, Computability is a special way of studying functions. To continue with our analysis we thus have to introduce a significant framework.

Formally in Computability Theory, when we speak of functions we are talking of a map whose domain is a subset of $\mathbb{N}^k$ for some $k \geq 1$ and whose range is $\mathbb{N}$. It is not assumed to be surjective or injective, but should be well-defined. It is rare for us to specify what the exact domain of these functions are, but we use the follow notation to specify where they are defined. For a function with input $\phi(x)$, we write $\downarrow$ to indicate that $x$ is in the domain and thus that this output is defined. When we later introduce Turing machines, we will often say that this arrow means the machine ``halts". On the other hand, if we write $\phi(x) \uparrow$, we mean that $x$ is not in the domain and the function is not defined on input $x$.

\begin{dfn}
A function is called ``partial recursive" if it is computed by a Turing machine\footnote{There are equivalent models of computation, such as via primitive recursion or unlimited register machines which can also be used as a definition, however we do not need to discuss them here.}. By computed, we mean that the function will return the same output for a given input as the machine. The precise definition of Turing machines will not be given, but can be found in \cite[p. 7-9]{soare2016turing} or any other introductory text book on the subject. 
\end{dfn}

\begin{dfn}
A partial recursive function is called ``total recursive", ``recursive" or ``computable" if it halts ($\downarrow$) on all inputs in $\mathbb{N}^k$ for the appropriate input arity $k$. 
\end{dfn}

We do not worry about the specific details of Turing machines, as we can take a higher level approach to keep results simple. The Church-Turing Thesis, which is near universally accepted and assumed here, states that Turing machines fully capture what is meant by an ``algorithm" informally. This is stated in \cite[p. 7]{soare2016turing}, under the name Turing's Thesis. Thus, if $f$ copies the behaviour of a set of fully specified steps to perform on inputs transforming them in finitely many ways before returning outputs, then a Turing machine exists to compute this $f$. Additionally $f$ can be permitted to have infinitely many steps, but then it must return no output on this input. This matches with our intuition that any process taking infinite time cannot also finish and give an output. As such we can consider ``computable" for functions as equivalent to saying that an informal total algorithm exists. \par
The only details we use about Turing machines is that they are finite in ``size", which allows them to be enumerated, and that they deal with finite inputs and outputs. There are several ways to perform this enumeration, but we take the canonical G\"{o}del numbering, which maps each machine to a natural number $e$ by some specific procedure which always halts. To read more about this procedure see \cite[p. 181]{openlogicproject2022sample}, although this may be hard to approach without first understanding its previous sections on Turing machines.

\begin{dfn}
If $e$ is the G\"{o}del number of a machine, we say that $\varphi_e$ is the partial recursive function corresponding to $e$. Additionally, Turing machines can be informally taken to accept multiple inputs such as via the same coding methods we described. As such when the machine accepts $n > 1$ inputs coded by $\langle v_1, \dots, v_n \rangle$, we specify the function $\varphi_e^{(n)}$ which takes multiple inputs $(v_1, \dots, v_n)$ rather than the single input $ \langle v_1, \dots, v_n \rangle$. We abuse notation by dropping the superscript indicating arity when it obvious because the inputs are specified to avoid cluttered formulae.
\end{dfn}

We can now see that the Church-Turing Thesis is equivalent to guaranteeing $\forall f$, a computable function, $\exists e \in \mathbb{N}$ s.t. $\varphi_e = f$. We can easily make our choice of $e$ unique, for example by just choosing the least $e$ which is a valid code, and this choosing is a computable procedure. This uniqueness is assumed without statement where convenient, so that we don't have to consider families of codes, but just a single representative.

\begin{dfn}
A set $A$ is called computable if its characteristic function $\chi_A$ is computable. We say that G\"{o}del number $e$ is the index of the set $A$ when $\varphi_e = \chi_A$. From a set's index, we can thus recover the set computably.
\end{dfn}

\begin{dfn}
For G\"{o}del number $e$, $W_e = \text{dom}\{ \varphi_e\} = \{ n \in \mathbb{N} : \varphi_e (n) \downarrow\}$ is the set of inputs with defined outputs. Then a set $A \subseteq \mathbb{N}$ is called ``recursively enumerable" (r.e.) if $A = W_e$. Clearly all computable sets are also r.e., as we can choose the function $\varphi_e = 1$ when $\chi_A = 1$ and doesn't halt otherwise, which is a computable procedure. Intuitively, if we given an $x \in A$ for a r.e. set, we can ``accept" the membership of $x$ in finite time. However, we cannot ``reject" some $y \notin A$ in finite time, unlike for computable sets, meaning we cannot be sure which elements belong to $A$.
\end{dfn}

 If the computable sets are those whose membership is determined by total computable functions, the r.e. sets are the domains on which partial machines halt. Intuitively we can think of a Turing machine which continuously prints its numbers, which are elements of the set it enumerates. If the set is r.e. but not computable (and therefore not finite), you can imagine this as an infinite process which never stops. Given time we may see that a certain element is in the r.e. set, but we can never know for sure that an element is \textit{not} in the set; it may just be that the machine hasn't gotten around to enumerating it yet.

 From here, we just establish a few classical theorems, taken without proof. These are all well-known and standard results, and are likely to be any introductory textbook on computability. We will reference \cite{soare2016turing} as a source for standardised notation.

\begin{thm}[Enumeration Theorem]\label{thm:enum}
For every G\"{o}del number $e$, there exists a partial recursive function $\varphi_z(e,x) = \varphi_e(x)$.
\end{thm}

We can interpret the Enumeration Theorem as positing the existence of a universal function (or universal Turing machine), which on input of a code $e$ and an input $x$ will simulate the run of the $e^\text{th}$ function on the input. It is easy to see that with our coding methods we can extend this to higher arities, but this suffices for our uses.

\begin{thm}[$s$-$m$-$n$ Theorem]\label{thm:smn}
For any choice of $m,n \geq 1$, there exists a bijective computable function $s^m_n$ of arity $m+1$ such that $\forall x,y_1, \dots, y_m, z_1, \dots, z_n,$
$$\varphi^{(n)}_{s^m_n(x, y_1, \dots, y_m)}( z_1  \dots, z_n) = \varphi^{(m+n)}_x(y_1, \dots, y_m, z_1  \dots, z_n).$$
\end{thm}

The notation of the \nameref{thm:smn} can be daunting, but intuitively it tells us that if we have a computable function of arity $m+n$, we can ``fix" some fraction $m$ of those inputs and compute the code of a new function with arity $n$ which computes the same function. Obviously the index of the new function will depend on what we fix our inputs to, which is why $s^m_n$ depends upon the fixed inputs as well as the code of the original function.

\begin{dfn}
A set $X$ is called an index set if $\forall x,y$ which are G\"{o}del numbers $x \in X \land \varphi_x = \varphi_y \implies y \in X$.
\end{dfn}

Here we are not assuming unique codes (indices), but returning back to coding of Turing machines where several machine can compute the same function. Therefore the $x$ and $y$ in the statement can be different. The importance of the definition is that index sets make claims about the ``features" of functions themselves, which thus means the index set includes all valid codes for the same function as they have the same features. Some typical examples might be the set of all codes for total functions, or the set of all codes for constant functions.

\begin{thm}[Rice's Theorem]\label{thm:Rice}
If $X \neq \emptyset$ and $X \neq \mathbb{N}$ is an index set, then $X$ is not computable.
\end{thm}

Put simply, Rice's Theorem tells us we cannot compute the index set for any non-trivial property. Looking at our previous examples, we see that at least some function is total, and at least some function isn't total, which now suffices to see that the index set of all total functions isn't computable. We use Rice's Theorem without proof, but one is given in \cite[p. 19]{soare2016turing}

\begin{dfn}\label{dfn:Reducible} 
A set $A$ is said to be many-one reducible to $B$ if there exists a computable $f$ such that $x \in A \iff f(x) \in B$. We denote this by writing $A \leq_m B$.
\end{dfn}

\begin{lem}\label{lem:Reduction}
If $A \leq_m B$ via function $f$, then $B$ is computable implies $A$ is computable. Further, $B$ is r.e. implies that $A$ is r.e.
\end{lem}
\begin{proof}
The first result is quite trivial; if $\chi_B = \varphi_e$, then $\chi_A(x) = \chi_B(f(x))$. The second result follows in much the same way. If $B = \text{dom}\{ \varphi_e\}$, then $A= \text{dom}\{ \varphi_e(f)\}$.
\end{proof}

The main function of reductions in computability theory is that they serve as a way to re-use our existing knowledge. Rather than showing fully that a set is computable, we can show it reduces to a computable set. We can also use the reduction to show non-computability. If a known non-computable set $K$ reduces to some other set $X$, we can be certain that $X$ is also non-computable.

\section{Second-Order Arithmetic}

Following the methodology of \cite{benedict2022arrow}, we can formalise the concepts introduced by \cite{kirman1972arrow} and \citeauthor{mihara1995arrow} within the framework of second-order arithmetic. As Eastaugh notes, second-order arithmetic is often used in mathematical logic and as such there is a significant number of existing techniques and results which can be applied here. Showing our computability-theoretic results hold in second-order arithmetic also enables further research in proof theory. For example, using the techniques of reverse mathematics, it could be possible to examine reversals of the main theorems of \citeauthor{mihara1995arrow}.

\begin{dfn}[Peano Axioms]\label{dfn:Peano}
Consider variables $x,y,z \in \mathbb{N}$ and symbols $(0,S,+,\times,<)$ which represent, in order: The zero symbol, the unary successor function $S(x)$, the binary addition function $+(x,y) \equiv x + y$, the binary multiplication function $\times(x,y) \equiv x \times y$ and the binary ``less than" relation $<(x,y) \equiv x < y$. The Peano arithmetic axioms minus induction (also called Robinsons's $Q$ axioms) require that
\begin{enumerate}[label=(Q\arabic*)]
    \item $\forall x,y, ( S(x) = S(y) \implies x = y )$ and
    \item $\forall x, \neg ( 0 = S(x) )$ and
    \item $\forall x, ( x = 0 \lor \exists y \text{ s.t. } ( x = S(y) ) )$ and
    \item $\forall x, ( x + 0 = x )$ and
    \item $\forall x,y, ( x + S(y) = S(x+y) )$ and
    \item $\forall x, (x \times 0 = x)$ and
    \item $\forall x,y, ( x \times S(y) = ( x \times y ) + x )$ and
    \item $\forall x,y, ( x < y \iff \exists z \text{ s.t. } ( s(z) + x ) = y )$.
\end{enumerate}
These arithmetic axioms aim to encapsulate basic mathematical intuition with respect to the natural numbers, and mean that the symbols given can be interpreted in the obvious way and have all expected properties.
\end{dfn}

\begin{dfn}
Second order arithmetic is a two-sorted language denoted $\mathcal{L}_2$. The first sort are the number variables $v_1, v_2, \dots$ whose domain is the set of natural numbers $\mathbb{N}$. The second sort are set variables $X_1, X_2, \dots$ whose domain is the power set of the natural numbers $\mathcal{P}(\mathbb{N})$. The symbols $(0,S,+,\times,<)$ can all be inherited from the \nameref{dfn:Peano}, as can the arithmetic axioms (minus induction) which dictate the behaviour of our operators w.r.t the number variables. Additionally we can define the set membership symbol $\in$ to express that a number variable belongs to a set variable. Finally we have quantifiers $\forall$ and $\exists$ which can quantify over number or set variables. These symbols alongside the logical symbols construct sentences over our language.
\end{dfn}

Up to this point, we have been defining our mathematical objects informally in terms of our intuitive understanding of objects like general sets, functions and structures in topology and algebra. In order to formalise these proofs in $\mathcal{L}_2$ for potential further analysis with respect to compatibility, we need to ensure our proofs hold in the basic theory and that these objects can be ``coded" into the domain of our number and set variables. Since the purpose of coding is to ``compress" otherwise high level mathematical objects into the two sorts provided by $\mathcal{L}_2$, it is important that our coding gives a unique representation, and that we can ``invert" our representation to recover the original object without loss. In set theory, there are several ways to represent a set as ``ordered" such that $\langle x,y \rangle \neq  \langle y,x \rangle$. The most commonly used is the Kuratowski definition that $\langle x,y \rangle = \{ \{ x,y \} , \{ x \} \}$.

\begin{dfn}
In arithmetic, a ordered pair $\langle x,y \rangle$ for $x,y \in \mathbb{N}$ can be coded by a single natural number $(x + y)^2 + x$. It is simple enough to see this coding is unique and invertible. An ``ordered n-tuple" or ``sequence" extends this definition by recursively\footnote{Here we are being informal in what is meant by recursion, but it is not really necessary. Instead there is a proof of the existence of the so-called beta function which codes sequences without requiring recursion. For more detail, see \cite[p. 350]{openlogicproject2022intermediate}} constructing $\langle v_1, v_2, \dots, v_n \rangle = \langle \langle v_1, \dots, v_{n-1} \rangle, v_n \rangle$.
\end{dfn}

\begin{dfn}
An n-ary relation $R$ can be interpreted as a set which contains $\langle v_1, \dots, v_n \rangle$ whenever $R(v_1, \dots, v_n)$ holds. This set of ordered tuples can then be coded by our previous definition to code $R$ as a set of natural numbers, meaning it corresponds to a set variable. A function $f:X \to Y$ is then just a binary relation s.t. if $\langle x,y \rangle \in f$, then $x \in X$ and $y \in Y$, which also satisfies $\forall x \in X, \exists ! y \in Y \text{ s.t. } \langle x,y \rangle \in f$. These functions are then set objects as they consist of a collection of codes for ordered pairs.
\end{dfn}

\section{Mihara's Result}

Unfortunately the machinery established so far doesn't allow us to fully encapsulate the results of \cite{kirman1972arrow} within $\mathcal{L}_2$. If sets of alternatives, $X$, and voters, $V$, are finite, we can represent them as tuples of natural numbers which can be coded to a single value. We can then consider functions over these sets meaning we can represent profiles and so on. \par
However, even lifting our restriction to the point of allowing only countable sets $V$ and $X$, we already see that profiles will be uncountable and so social welfare functions would not be second-order objects. \par
We address this issue in a few ways, following the methods of \cite{mihara1997arrow}. The next section first introduces the concept of algebras, which are a useful tool for limiting our domain. We are then able to begin to explore our main results by establishing our modified conception of what it means to be computable.

\begin{dfn}
A set $\mathscr{B}$ is called a boolean algebra\footnote{Boolean algebras can be defined in more general terms as structures equipped with several operations, but for our purposes this is sufficient since we are always working over the power set of $V$. It is not hard to see that our definition of boolean algebra meets the required conditions if we choose  operations to mimic the behaviour of set intersection, union and complementation.} over set $V$ if its elements are subsets of $V$ and $\forall A,B \in \mathscr{B}$ it is
\begin{enumerate}[label=(\roman*)]
    \item\label{alg:union} closed under union, meaning $A \cup B \in \mathscr{B}$, and
    \item\label{alg:intersect} closed under intersection, meaning $A \cap B \in \mathscr{B}$, and
    \item\label{alg:complement} closed under complementation, meaning $A^c \in \mathscr{B}$, and
    \item\label{alg:nonempty} $\emptyset, V \in \mathscr{B}$.
\end{enumerate}
\end{dfn}

\begin{dfn}
If $\mathscr{B}$ is a boolean algebra over $V$, we say that a set $A$ is ``$\mathscr{B}$-measurable" whenever $A \in \mathscr{B}$. Further we say a voter profile $f \in F$ is $\mathscr{B}$-measurable if and only if $\forall a,b \in X, \{ v \in V : af(v)b \} \in \mathscr{B}$.
\end{dfn}

\begin{dfn}
 We denote the set of all $\mathscr{B}$-measurable profiles by $F_\mathscr{B}$, and say a social welfare function $\sigma$ is a $\mathscr{B}$-social welfare function if its domain is $F_\mathscr{B}$.
\end{dfn}

We begin by introducing Mihara's main theorem. For the following section we fix the following: $X$ is a set of at least $3$ alternatives (countable or uncountable), $V$ is a set of voters which are represented by elements in $\mathbb{N}$, which can be some sort of unique identifier, such as a national insurance number in the UK. Clearly in most cases the voters will be a finite set, but as shown in \citeauthor{kirman1972arrow}, the finite case is uninteresting as dictatorship is already guaranteed, so we instead extend to the entirety of $\mathbb{N}$. \par
The boolean algebra over $\mathbb{N}$ that we choose is $\text{REC} = \{N \subseteq \mathbb{N} : N \text{ is a computable set} \}$, and it is not hard to see that this satisfies our requirements\footnote{ We can see if $\chi_A$ and $\chi_B$ are computable, then so is $\chi_{A \cap B} = \chi_A \times \chi_B$, $\chi_{A \cup B} = \chi_A + \chi_B - \chi_{A \cap B}$ and $\chi_{A^c} = 1 - \chi_A$, which gives the required closures.}. Thus our profiles will be REC-measurable if the set of voters (often referred to as a \textit{coalition} in the literature) which prefer any fixed alternative $x$ to another $y$ is itself computable. Correspondingly, the REC-social welfare functions are those which only accept REC-measurable inputs. This means that our social welfare function may be undefined for non-measurable inputs. \par
In order to consider social preferences in the setting of computability, it is clear we must have some restrictions placed on the ways in which voters can rank alternatives. The first is already given by our choice of REC as our algebra, as not all social welfare functions are REC-social welfare functions. Mihara justifies this restriction as reasonable, since requiring REC-measurability is equivalent to demanding that our voters are able to give an algorithm to describe their voting behaviour. As such we are only excluding indescribable voting behaviour, which seems natural in any applied setting. \par
This is not all; if our set $X$ of alternatives is not finite, it may not be enough of a restriction. As Mihara (1995) notes, $F_{\text{REC}}$, which is the set of all REC-measurable profiles, is itself uncountable:

\begin{lem}
$F_{\text{REC}}$ is uncountable if $X$ is infinite.
\end{lem}
\begin{proof}
Our proof hinges on showing there is a bijection from $\mathcal{P}(X)$, which is uncountable since $X$ is infinite, to the set of preferences. For a subset $A \in \mathcal{P}(X)$, the preference is constructed s.t. $xPy$ iff $x \in A$ and $x \notin A$, which is clearly one-to-one. It is immediately clear that this satisfies \ref{dfn:weakorderasym}. Further we assume $\neg xPy \land \neg yPz$ and for a contradiction $xPz$. So $x \in A$ and $z \notin A$; now we can be sure that $y \in A$, else we get $xPy$, but then since $y$ and $z$ lie in and out of $A$ respectively, we get a contradiction since $yPz$. Thus we satisfy \ref{dfn:weakordernegtrans}, and so our construction is really a preference. Now the profile in which all voters have the same preference is clearly REC-measurable, as $\{ v \in V : af(v)b \}$ is either the entire set $\mathbb{N}$ or $\emptyset$, which are in REC by algebra property \ref{alg:nonempty}. This serves as an injection from preferences to profiles. Thus this forms an uncountable subset of $F_{\text{REC}}$, which is therefore itself also uncountable.
\end{proof}

This uncountability result is immediately problematic. It is impossible to talk about the computability of social welfare functions if they have an uncountable domain, as we have no way to code the all inputs of the function as finite inputs for algorithms, which are inherently countable; they cannot be said to correspond to either computable or non-computable sets without further considerations. The importance of this problem cannot be overstated; whilst \citeauthor{mihara1995arrow} is relatively nonchalant about things, this is the primary motivating factor for needing an alternative definition for computability in his paper. Mihara's proposed  solution to this isn't simply to restrict to a countable domain, so we first define how to code our profiles, then define our modified notion of computability. \par
From this point onward, we proceed with the proof by assuming the first four of \nameref{dfn:arrow}, \ref{A1:alternatives}, \ref{A2:welfarefunction}, \ref{A3:unanimity} and \ref{A4:independence}. This is generally because our proof strategy will be to assume these as a premise and aim for a contradiction with \ref{A5:nondictatorial}, showing dictatorship. In particular, we use the axiom of independence, \ref{A4:independence}. A corollary of independence as defined is that in order to determine the behaviour of any social welfare function $\sigma$ on distinct alternatives $\{x,y\}$ for profile $f$, we only need to know the behaviour of $f$ restricted to $\{x,y\}$. This gives us a strategy for understanding the behaviours of REC-social welfare functions restricted to pairs by looking at codes for profiles.

\begin{dfn}\label{dfn:profilecode}
For a pair of distinct elements $x,y \in X$, the REC-Measurable profile $f \in F_\text{REC}$ restricted to $\{x,y\}$ is coded by three terms; $e_1$ indexes the set $\{v \in V : x f(v) y \}$, $e_2$ indexes the set $\{v \in V : y f(v) x \}$ and $e_3$ indexes the set $\{v \in V : \neg (x f(v) y \lor y f(v) x) \}$. We simplify this triple to the single code $e = \langle e_1, e_2, e_3 \rangle$.
\end{dfn}

\begin{lem}
The code in Definition \ref{dfn:profilecode} is well-defined. 
\end{lem}
\begin{proof}
Clearly $e_1$ and $e_2$ exist since REC-measurability ensures the sets they index are in REC, and thus are computable. We can thus be sure that $e_3$ exists, as by construction it is just the complement of the union of the two computable sets indexed by $e_1$ and $e_2$. This is then in REC by the algebra properties \ref{alg:union} and \ref{alg:complement}.
\end{proof}

It is not necessary to have all three codes, as it is not hard to see that this is a partition, and thus any third set can be computed from the other two. However, we adopt this extra notation to follow in the conventions of Mihara.

\begin{dfn}[Pairwise Computable]\label{dfn:PC}
A (REC-)social welfare function $\sigma$ is called ``pairwise computable" or ``PC" if, for each distinct ordered pair $\langle x,y \rangle \in X \times X$, there exists $\gamma_{\langle x,y \rangle}$ a partial recursive function and for each REC-measurable profile $f \in F_\text{REC}$ coded by $e = \langle e_1, e_2, e_3 \rangle$ on $\{x,y\}$, we have that $\gamma_{\langle x,y \rangle}$ satisfies
\begin{align*}
x \sigma(f) y       & \implies \gamma_{\langle x,y \rangle} (e) = 1, \text{ and } \\
\neg x \sigma(f) y  & \implies \gamma_{\langle x,y \rangle} (e) = 0.    
\end{align*}
This means that on input of a valid code $e$ for a profile on $\{x,y\}$, our partial recursive function $\gamma$ will tell us the relative position of the two alternatives in the welfare function's verdict. This means that even if we cannot compute the entire social welfare function, we can still query whether any alternative is preferred over another.\par
It is of note that each $\gamma$ only has its behaviour specified on the input of valid codes for profiles, and may return anything for other natural numbers. This isn't much of a restriction, as any correctly organised vote will produce a valid code. So whilst problems can be considered in theory, they cannot exist in practice except by human error; we have yet to find a cure for human fallibility, computable or not.
\end{dfn}

\begin{thm}[Mihara's Theorem 1]\label{thm:mihara}
Let $\sigma : F_\text{REC} \to W$ be a REC-social welfare function. Assume the first four of \nameref{dfn:arrow}, \ref{A1:alternatives} - \ref{A4:independence}. Then if $\sigma$ is \nameref{dfn:PC}, then $\sigma$ satisfies $\neg$\ref{A5:nondictatorial} and is dictatorial.
\end{thm}

Before we can prove this theorem, we develop considerably machinery to help us, through a series of lemmas.

\begin{lem}\label{lem:mihara1}
There is a bijective computable function $r$ such that for all G\"{o}del numbers $e$ and for all $n \in \mathbb{N}$
$$
\varphi_{r(e)}(n) = 
\begin{cases}
1           & \text{if } \varphi_e(n) = 0, \\
0           & \text{if } \varphi_e(n) \downarrow \text{ and } \varphi_e(n) \neq 0, \\
\uparrow    & \text{if } \varphi_e(n) \uparrow.
\end{cases}
$$
Further if $e$ is the index of set $A \subseteq \mathbb{N}$, $r(e)$ is the index of $A^c$ relative to $\mathbb{N}$.
\end{lem}
\begin{proof}
We define truncated subtraction given by the symbol $\dotdiv$ (also called monus) as subtraction for natural numbers in the obvious way, but if $a - b < 0$, we instead let $a \dotdiv b = 0$. It is well known that this operation is computable, and is a good first exercise for any reader more interested in exploring the details of computability. We can see that if we let $\psi(e,n) = 1 \dotdiv \varphi_e(n)$, the \nameref{thm:enum} ensures $\psi$ is partial recursive, with code $z$. Then by using the \nameref{thm:smn} we can see that there exists a bijective functions $r$ as desired s.t. $\varphi_{r(z,e)}(n) = \varphi_z(e,n) = 1 \dotdiv \varphi_e(n)$. To spell the use of this theorem out, we have $\varphi^{(1)}_{s^1_1(z,e)}(n) = \varphi^{(2)}_z(e,n)$ for $m = n = 1$, $r := s^1_1$, $x := z$, $y_1 := e$ and $z_1 := n$ if we mimic the notation of the statement explicitly. \par
By ignoring the input $z$ we can treat $r$ as a function of just $e$. It is simple to check by case analysis this has the desired property. We can see that if $e$ indexes $A$, then $\varphi_e = \chi_A$, and so $\varphi_{r(e)} = 1 - \chi_A = \chi_{A^c}$, so $r(e)$ indexes the complement.
\end{proof}

\begin{lem}\label{lem:mihara2}
We define the set of indices of computable sets by $\text{Rec} = \{ e : \exists A \in \text{REC} \text{ indexed by } e \}$\footnote{This notation mimics that of \cite[p. 17]{soare2016turing}, rather than \cite[p. 271]{mihara1997arrow} who calls it CREC. Either way, it is not to be confused with REC which contains the computable sets themselves rather than Rec which contains the indices}. Then Rec is not r.e.
\end{lem}
\begin{proof}
Fix a $\Sigma_2$ set\footnote{The notation used in this proof will not be introduced until later at \ref{dfn:hierarchy}, because it requires a deeper understanding of computability theory. For now understand this means $\exists$ a computable relation $R$ and $x \in A \iff \exists y_1 \forall y_2 R(x,y_1,y_2)$ holds as a relation.} $A$. By \cite[p. 66\footnote{Unfortunately I could not access this text, nor find the original version cited by \cite[p. 271]{mihara1997arrow}; as such the page numbering may have changed in newer editions of \cite{soare1999recursively}.}]{soare1999recursively}, there exists a computable $f$ s.t.
\begin{align*}
    e \in A \implies    & \varphi_{f(e)}(n) \downarrow \text{ for only finitely many inputs $n \in \mathbb{N}$ and } \\
    e \notin A \implies & \varphi_{f(e)}(n) = 0, \text{ which is clearly computable}.
\end{align*}
We now assume for contradiction that Rec is r.e., meaning that there is a G\"{o}del number $\overline{e}$ s.t. Rec $ = \text{dom}(\varphi_{\overline{e}} )$, and $\varphi_{\overline{e}}$ halts whenever our input is in Rec. Given any code $e$, if $e \notin A$, $f(e)$ is returning the index of a computable set, so $f(e) \in \text{Rec}$. On the other hand, if $e \in A$, $f(e)$ is returning the index of a partial but not computable set (since it must fail to halt after finitely many inputs), so $f(e) \notin \text{Rec}$. \par
Therefore $e$ is in $A^c$ if and only if $f(e)$ is in Rec, and thus $A^c$ is many-one reducible to Rec, written $A^c \leq_m \text{Rec}$. For contradiction, we assume Rec is r.e., which is equivalent\footnote{We return to make this in \ref{lem:postHierarchy}.} to assuming it is in $\Sigma_1$. Then by properties of the reduction in \ref{lem:reduceHierarchy}, Rec $ \in \Sigma_1 \implies A^c \in \Sigma_1 \implies A \in \Pi_1$. But it is known that there exist sets strictly in $\Sigma_2$ not in $\Pi_1$\footnote{This comes about from the properness of the hierarchy, established in \cite[p. 85]{soare2016turing}.}, which gives a contradiction, and Rec is not r.e.
\end{proof}

\begin{lem}\label{lem:mihara3}
For a given distinct pair $x,y \in X$, we define the set of codes of REC-measurable profiles restricted to $\{x,y\}$ by
$$S_{\{x,y\}} = \{ e : \exists f \in F_\text{REC} \text{ coded by } e \text{ restricted to } \{x,y\} \}.$$ Then $S_{\{x,y\}}$ is not r.e.
\end{lem}
\begin{proof}
We begin by fixing some $x,y$, then assuming for contradiction that $S_{\{x,y\}} = S$ is r.e. Now we define the computable function $r$ as in \ref{lem:mihara1}, and fix $e_3$ as an index for the empty set (which is clearly computable). Consider $e_1 \in \text{Rec}$, and the triple $e = \langle e_1, r(e_1), e_3 \rangle$; clearly if $e \in S$, since our coding process is invertible by a simple computable method we can recover $e_1$ which is guaranteed to index the set $\{v \in V : x f(v) y \} \in REC$. Thus $e_1 \in Rec$. \par 
On the other hand, if we assume $e_1 \in Rec$ and indexes the computable set $A$, then we define a profile $f$ restricted to $\{x,y\}$ s.t. $A = \{v \in V : x f(v) y \}$ and $A^c = \{v \in V : y f(v) x \}$, both of which are REC-measurable. Furthermore, by \ref{lem:mihara1}, we know $r(e_1)$ is the index of $A^c$. This profile is indifferent on no elements, so $\{v \in V : \neg (x f(v) y \lor y f(v) x) \} = \emptyset$ which is indexed by $e_3$. Since all the sets are computable, $f \in F_\text{REC}$ and has code $e$, and so $e \in S$. Therefore $e_1 \in \text{Rec} \iff e = \langle e_1, r(e_1), e_3 \rangle \in S$. \par
As noted, constructing $e$ from $e_1$ is computable by some function, meaning Rec is many-one reducible to $S$, and so $\text{Rec} \leq_m S$. Now since $S$ is assumed to be r.e., by properties of the reduction given in \ref{lem:Reduction}, we get that Rec is r.e. But this is of course in contradiction to \ref{lem:mihara2}. Therefore $S_{\{x,y\}}$ is not r.e.
\end{proof}

From here we shall re-use some of the results from the \citeauthor{kirman1972arrow} proof. In his original paper, \citeauthor{mihara1995arrow} takes these results from \citeauthor{armstrong1980arrow}, who himself says they are essentially taken from \citeauthor{kirman1972arrow}. However \citeauthor{kirman1972arrow}'s original proof is likely to be more readable as their notation and setup is similar to that used here. \par
The only real modification is that \cite{armstrong1980arrow} talks about ultrafilters over an algebra, which is what we shall apply here. This is necessary as we can't generally code all filters without restricting them to our computable domain of the algebra.

\begin{dfn}
Given a filter $\mathcal{F}$ over a set $V$ and a boolean algebra $\mathscr{B}$ on the same set, $\mathcal{F}$ is said to be a filter on $\mathscr{B}$ when $\mathcal{F} \subseteq \mathscr{B}$. A filter is then an ultrafilter over the algebra if it is maximally strict with respect to filters over the algebra, just as in the ordinary case.
\end{dfn}

This is a very natural extension since any filter over an algebra is a filter in the ordinary sense, and a filter in the ordinary sense is a filter over the algebra if we take it's intersection with the algebra. We only need to take care with regard to ultrafilters, as maximality w.r.t. an algebra may not be maximal in the ordinary sense; the other direction does however hold. A previously maximal filter will still be maximal in the algebra. This means that we can use the general statements by Kirman and Sondermann in \ref{thm:ks1} and \ref{thm:ks2} modified for ultrafilters on a boolean algebra without re-proof. It also means all our expected properties of ultrafilters still hold.

\begin{dfn}
Given $\sigma : F_\text{REC} \to W$, a REC-social welfare function which also satisfies \ref{A1:alternatives} - \ref{A4:independence}, we let $g(\sigma) = \mathcal{U}_\sigma$ be an ultrafilter as defined in \ref{thm:ks1}. We then define the $\beta_\sigma$ partial function on $\mathbb{N}$ s.t.
$$
\beta_\sigma (e) = 
\begin{cases}
1           & \text{if $e$ is the index for a recursive set in $\mathcal{U}_\sigma$} \\
0           & \text{if $e$ is the index for a recursive set not in $\mathcal{U}_\sigma$} \\
\uparrow    & \text{if $e$ is not the index for a recursive set.} \\
\end{cases}
$$
Since indices represent at most one set under our G\"{o}del numbering, we don't have to worry about the function being poorly defined. At this point it is unclear as to whether our function is partial recursive or not.
\end{dfn}

\begin{dfn}[Decidability of Decisive Coalitions]\label{dfn:DDC}
We say $\sigma$ satisfies ``Decidability of Decisive Coalitions" (DDC) if $\beta_\sigma$ as defined has an extension to a partial recursive function.
\end{dfn}

\begin{lem}\label{lem:mihara4}
Given $\sigma : F_\text{REC} \to W$, a REC-social welfare function which also satisfies \ref{A1:alternatives} - \ref{A4:independence}, then $\sigma$ satisfies $\neg$\ref{A5:nondictatorial} and is dictatorial if and only if it satisfies DDC.
\end{lem}
\begin{proof}
First suppose that $\sigma$ satisfies $\neg$\ref{A5:nondictatorial}. Then by \ref{thm:ks2}, the ultrafilter $\mathcal{U}_\sigma$ on REC is fixed. Examining the proof of \ref{lem:singletongenerated} a little, we can remark that in the ``only if" direction we only use the fact that the ultrafilter is finite to show it is fixed. As such we can use the same method to prove all fixed ultrafilters are of the form $\{ U \subseteq V : v_0 \in U\}$, with no assumption that $V$ is finite. Now we know $\exists v_0 \in V$ s.t. $\mathcal{U}_\sigma = \{ U \in REC : v_0 \in U\}$. Since $v_0$ is now our dictator element, we can query the verdict of the social welfare function by just checking the dictator's choice. Thus we define
$$
\overline{\beta_\sigma} (e) = 
\begin{cases}
\varphi_e(v_0)  & \text{if $\varphi_e(v_0) = 0$ or $1$,} \\
\uparrow        & \text{otherwise.} \\
\end{cases}
$$
By the \nameref{thm:enum}, we know that $\overline{\beta_\sigma}$ is partial recursive since $\varphi_e$ is always partial recursive. It then suffices to observe that $\overline{\beta_\sigma} = \beta_\sigma$, since whenever $e$ indexes a recursive set $A$, it is enough to query $\varphi_e(v_0) = \chi_A(v_0)$, which will return $1$ if $v_0 \in A \implies A \in \mathcal{U}_\sigma$ and $0$ if $v_0 \notin A \implies A \notin \mathcal{U}_\sigma$. So the DDC property holds. \par
In the other direction, we begin by assuming the DDC property. For a contradiction, we assume \ref{A5:nondictatorial}, i.e. that the REC-social welfare function is non-dictatorial. Therefore by \ref{thm:ks2}, we know $\mathcal{U}_\sigma$ is a free ultrafilter. Now we can see that any ultrafilter which contains any finite set must be fixed, because otherwise we could add any singleton from the finite set and apply closure to get a finer filter. Therefore the contrapositive states that any free filter contains no finite sets. But since it contains every set or its complement, it must contain all cofinite sets, which are those whose complements are finite. \par
At this point we take what seems to be a slight tangent, but which returns to complete the proof: First define $\text{K} = \{ e : e \in W_e\}$, often called G\"{o}del's diagonal set. Unfortunately this is not an index set, so we cannot apply Rice's Theorem. Instead, we can show it is non-computable by more elementary means of diagonalisation, such as in \cite[p. 13-14]{soare2016turing}, although we are more explicit. If K were to be computable, we could simply construct a new function;
$$
g (e) = 
\begin{cases}
\varphi_e(e) + 1  & \text{if $e \in $ K,} \\
0        & \text{if $e \notin $ K.} \\
\end{cases}
$$
Now $g$ is computable using the indicator of K, so this means $g = \varphi_k$ for some $k \in \mathbb{N}$. But then if we consider along the ``diagonal", we find a problem at $g(k)$. If $k \in K$, $g(k) = \varphi_k(k) + 1 = \varphi_k(k)$, which is impossible. If $k \notin K$, $g(k) = 0$, but this implies $k \in \text{dom}\{ \varphi_k\} = W_k$ which then implies $k \in K$ for a contradiction. So we can see K is not computable. \par
We can then see it is r.e. because $e \in \text{K}$ if and only if $\exists x \exists y$ s.t. $\varphi_x(e)=y$, which is a $\Sigma_0^1$ statement. Thanks to \cite[p. 20\footnote{As before, I could not access the original. For the initial citation see \cite[p. 273]{mihara1997arrow}.}]{soare1999recursively}, there is therefore a computable set (relation) $R \subseteq \mathbb{N} \times \mathbb{N}$ s.t. $e \in \text{K}$ if and only if $\exists z$ s.t. $R(e,z)$ holds. We then define for a fixed $z$ the relation
$$
h(e,u) = 
\begin{cases}
1   & \text{if $\exists z \leq u$ s.t. $R(e,z)$,} \\
0   & \text{otherwise,} \\
\end{cases}
$$
which is clearly computable since the quantifier is bounded (and thus a finite ``search"). Using the \nameref{thm:smn}, we have a computable $f(e) = s(z,e)$ s.t. $\varphi_{s(z,e)}(u) = \varphi^{(1)}_{s^1_1(z,e)}(u) = \varphi_z^{(2)}(e,u) = h(e,u)$. Now $e \in \text{K}$ implies some $z_0$ satisfies $R(e,z_0)$. So $\varphi_{f(e)}(u) = h(e,u) = 0$ only on finitely many inputs when $u < z_0$.
Therefore $f(e)$ is the index for a cofinite set. On the other hand, if $e \notin \text{K}$, $\varphi_{f(e)}(u) = h(e,u) = 0$ always, and so $f(e)$ is the index for $\emptyset$. Finally, using the DDC property, we can compute $\beta_\sigma(f(e))$ by composition. But $e \in \text{K}$ now implies $f(e)$ is cofinite which implies $\beta_\sigma(f(e)) = 1$. On the other hand, $e \notin \text{K}$ now implies $f(e)$ is the $\emptyset$ index which implies $\beta_\sigma(f(e)) = 0$. So interpreting inputs as codes for indices, $\beta_\sigma(f(e))$ is actually the indicator of the set K. Therefore K is computable, which is an obvious contradiction. Therefore we have $\neg$ \ref{A5:nondictatorial}, so the function has a dictator.
\end{proof}

Having shown that the DDC property is equivalent in this context to being dictatorial, we can finally complete our proof of \nameref{thm:mihara}. We do so simply by showing that REC-social welfare functions that are \nameref{dfn:PC} satisfy DDC.

\begin{proof}[Proof of \ref{thm:mihara}]
We assume the REC-social welfare function $\sigma : F_{\text{REC}} \to W $ satisfies \ref{A1:alternatives} - \ref{A4:independence}, and also is \nameref{dfn:PC}. For a fixed pair $\{x,y\}$ in our set of alternatives $X$, since the function is PC we have another partial recursive function $\gamma_{\langle x,y \rangle}$ which, on valid input $e$ coding a REC-measurable profile $f$, will return $1$ when $x \sigma(f) y$ and $0$ when $\neg x \sigma(f) y$. We then fix an index $e_3$ for $\emptyset$. By \ref{lem:mihara1}, we know there exists a bijective computable function $r$ with properties described. Then the triple $e = \langle e_1, r(e_1), r_3 \rangle$ is a valid profile (which is equivalent to being in $S_{\{x,y\}}$) whenever $e_1 \in \text{Rec}$, which was proven as part of the proof of \ref{lem:mihara3}. \par
Now if $e_1 \in $ Rec, $\beta_\sigma(e_1) = 0$ or $1$ depending on whether it is in the ultrafilter $\mathcal{U}_\sigma$. We let the computable set indexed by $e_1$ be called $A$, and by the results of \ref{lem:mihara1} we know $r(e_1)$ indexes the set $A^c$. Here we break into cases:
\begin{enumerate}
    \item Suppose $\beta_\sigma(e_1) = 1$. Then by construction, we see $A \in \mathcal{U}_\sigma$; looking back at our definition of the triple coding in \ref{dfn:profilecode}, we see that $A = \{v \in V : x f(v) y \}$. Now by \ref{thm:ks1}, since $A$ is in the ultrafilter, it is a sufficient coalition to decide the verdict, as $x f(A) y \implies x \sigma(f) y$. Finally, by the \nameref{dfn:PC} property, we see $\gamma(e) = 1$.
    \item Suppose $\beta_\sigma(e_1) = 0$. Similarly $A \notin \mathcal{U}_\sigma \implies A^c \in \mathcal{U}_\sigma$ by \ref{lem:ultracomplement}, and $A^c = \{v \in V : y f(v) a \}$. Again it is a sufficient coalition to decide the verdict, so $y f(A^c) x \implies y \sigma(f) x \implies \neg x \sigma(f) y$ by the asymmetry property \ref{dfn:weakorderasym}. Finally, by the \nameref{dfn:PC} property, we see $\gamma(e) = 0$.
\end{enumerate}
We've thus seen for inputs of $e_1 \in$ Rec, $\beta_\sigma(e_1) = \gamma(e)$. Since we can computably decode $e$ to re-obtain $e_1$, the fact that $\gamma$ is partial recursive implies that an extension of $\beta_\sigma$ is also partial recursive. But this is exactly the DDC criterion. So by \ref{lem:mihara4}, $\sigma$ is dictatorial.
\end{proof}

\section{Beyond the Results}

At this point, we have covered what we aimed to achieve in the paper, and have established the main result of \citeauthor{mihara1995arrow}. It is worth recapping the overall project. We initially established the impossibility and possibility results. We then introduced formal logic as a framework and computability as a methodology for analysis. We then used the definition of PC to show that no \nameref{dfn:PC} social welfare function is non-dictatorial. \par
Returning to the initial project of \citeauthor{fishburn1970arrow}, we have to consider the impact of our results. If all PC social welfare functions are dictatorial, the only non-dictatorial social welfare functions over REC as posited by \nameref{thm:fishburn} would be non-PC. Looking at the negation of \nameref{dfn:PC}, this means for any pair of candidates $\langle x,y \rangle$ we have no recursive method for determining their relative placements in the verdict. If we think about this, it means that we really have no decidable means to actually access the verdict, which seems to undermine the purpose of the vote entirely. After all, what is point of voting and coming to a decision if we have no computable means of agreeing on the outcome. Even if some human, non-computable meaning could be read from the social welfare function, it would be wrought with disagreements and claims of dishonesty. After all, who would trust a governing body that cannot prove why a certain candidate beat another. \par

Whilst this wraps up most essential discussion, we have only described computability in terms of the result of computations of Turing machines. This is natural, because as stated earlier the Church-Turing Thesis posits that this is the limit of ordinary computation. In the study of computability, however, it is very common to run into these non-computable sets. It is even common to see sets that are neither r.e. nor co-r.e. either. In this case, it is natural to ask whether we can quantify in some way the extent to which certain sets are non-computable. This motivates the notion of relative computability. Here we ask ourselves what we can compute if we just ``pretend" some set is computable by an unexplained means, which is commonly referred to as an ``oracle". We formalise this concept with the following definition.

\begin{dfn}
Given a set $X \subseteq \mathbb{N}$ which we call an oracle, a function is called ``$X$-partial recursive" if it is computed by a Turing machine which additionally during its run can query whether an element $n \in X$ and receive a $1$ or $0$ accordingly (and correctly). $X$-computable functions are just $X$-partial functions halting on all inputs.
\end{dfn}

It is typical of course for our choice of $X$ to be some non-computable set. After all, if we choose a computable set $X$, we don't change what is possible for our machine, as it could always have queried membership of $X$ by just computing it before. However, given non-computable sets as oracles, we can compute strictly more sets. At a minimum we can very uncreatively compute $X$ itself. \par
Since these machines are constructed in the same way as before, with the single additional property of querying membership of $X$, it is not hard to modify our G\"{o}del numbering to code these machines instead. We shall not consider the details of this process as this would require an explanation of G\"{o}del numbering to begin with. Instead, we assume for each $X$ there is another canonical coding for the oracle Turing machines.

\begin{dfn}
If $e$ is the G\"{o}del number of an $X$-oracle machine, we say that $\varphi_e^X$ is the $X$-partial recursive function corresponding to $e$. A set $A \subseteq \mathbb{N}$ is an $X$-computable set if it is the range of an $X$-computable function.
\end{dfn}

\begin{dfn}
For G\"{o}del number $e$, $W_e^X = \text{dom}\{ \varphi_e^X\} = \{ n \in \mathbb{N} : \varphi_e^X (n) \downarrow\}$. Then a set $A \subseteq \mathbb{N}$ is called ``$X$-recursively enumerable" ($X$-r.e.) if $A = W_e^X$.
\end{dfn}

\begin{dfn}
The ``Turing jump" of $X$, denoted by $X'$, is the set $\text{IND}^X = \{ e : e \in W_e^X\}$. We define the $n^\text{th}$ jump, denoted by $X^{(n)}$, via recursively applying the jump operation $n$ times so that $X^{(n+1)} = (X^{(n)})'$.
\end{dfn}

As we remarked earlier, if $X$ is chosen to be computable our machines are no more powerful than before. As such, we often consider all our ordinary definitions of computability as valid for $\emptyset$ as a canonical oracle. As such $\emptyset^{(n)}$ is taken to establish a natural hierarchy of sets. The reason some of these concepts are useful is because they tie back to another hierarchy which we made mention of earlier.

\begin{dfn}\label{dfn:hierarchy}
A set is in the classes $\Sigma_0, \Pi_0$ and $\Delta_0$ if it is computable. For $n>1$, we say a set $A$ in a class if there exists a computable relation $R$ such that $x \in A \iff \dots \exists y_k \forall y_{k+1} \exists y_{k+1} \forall y_{k+2} \dots R(x,y_1, \dots y_n)$ holds, where the $\forall$ and $\exists$ symbols are alternating, and there are $n$ of them. If the alternating sequence begins with the universal $\forall$ quantifier, it is in $\Pi_n$, and if it begins with the existential $\exists$ quantifier, it is in $\Sigma_n$. Then $\Delta_n = \Sigma_n \cap \Pi_n$.
\end{dfn}

We can make a few remarks on this hierarchy. First, we know through our knowledge of coding that the relation $\exists y_0 \exists y_1 R(y_0,y_1)$ can be coded to hold when another relation $\exists y S(y)$ holds where $y$ is a function of $y_0, y_1$. So, we can imagine that in our alternating sequence, we could have multiple existential or universal quantifiers in a row. Next, we can see that there is a relation between the $\Sigma_n$ and $\Pi_n$ sets. By taking the complement of one, we reach the other, as the negation of an existential quantifier is universal, and vice versa. The next remark is a very powerful result, thanks to Post and given in \cite[p. 84]{soare2016turing}. 

\begin{lem}\label{lem:postHierarchy}
$\emptyset^{(n)}$ relates to $\Sigma_{n+1}$ in such a way that the $\Sigma_{n+1}$ set are precisely those that are $\emptyset^{(n)}$-r.e. This is why in \ref{lem:mihara2} we say that the $\Sigma_1$ sets are exactly the r.e. sets.
\end{lem}

Finally, in \cite[p. 81]{soare2016turing} we also see that the previously defined reduction in \ref{dfn:Reducible} can used to interact with the hierarchy. Thus we have the following lemma.

\begin{lem}\label{lem:reduceHierarchy}
If $B \leq_m A$ and $A \in \Sigma_n$, then $B \in \Sigma _n$.
\end{lem}

Having established these final results, there is nothing left unexplained, and the dissatisfied reader can return to \ref{lem:mihara2} to convince themselves of the results. \par
However, there is still room for further exploration of these topics. Having introduced the notion of relative computability, it may be natural to ask whether we can use this tool of abstraction to reach a more general form of \nameref{thm:mihara}, in which we do not have REC-measurable sets, but sets measurable to an oracle $X$, and where $X$-REC can be defined as the set of $X$-computable functions. This has not been attempted in any of the surveyed literature, but we might expect that it should still hold. This is in part motivated by the fact that in \cite[p. 88]{soare2016turing}, we already see Rec is $\Sigma_3$ complete, and so $X$-Rec should be $\Sigma^X_3$ computable. This gives us hope that this abstraction might well be possible for any keen reader looking for further research.

\section{Bibliography}

\bibliographystyle{plainnat}
\nocite{*}
\bibliography{bib.bib}

\begin{thebibliography}{14}
\providecommand{\natexlab}[1]{#1}
\providecommand{\url}[1]{\texttt{#1}}
\expandafter\ifx\csname urlstyle\endcsname\relax
  \providecommand{\doi}[1]{doi: #1}\else
  \providecommand{\doi}{doi: \begingroup \urlstyle{rm}\Url}\fi

\bibitem[Armstrong(1980)]{armstrong1980arrow}
Thomas~E Armstrong.
\newblock Arrow's theorem with restricted coalition algebras.
\newblock \emph{Journal of Mathematical Economics}, 7\penalty0 (1):\penalty0
  55--75, 1980.

\bibitem[Armstrong(1985)]{armstrong1985precisely}
Thomas~E Armstrong.
\newblock Precisely dictatorial social welfare functions: Erratum and addendum
  to ‘arrows theorem with restricted coalition algebras’.
\newblock \emph{Journal of Mathematical Economics}, 14\penalty0 (1):\penalty0
  57--59, 1985.

\bibitem[Eastaugh(2022)]{benedict2022arrow}
Benedict Eastaugh.
\newblock Arrow to infinity: the reverse mathematics of social choice theory.
\newblock unpublished, 2022.

\bibitem[Fishburn(1970)]{fishburn1970arrow}
Peter~C Fishburn.
\newblock Arrow's impossibility theorem: concise proof and infinite voters.
\newblock \emph{Journal of Economic Theory}, 2\penalty0 (1):\penalty0 103--106,
  1970.

\bibitem[Kirman and Sondermann(1972)]{kirman1972arrow}
Alan~P Kirman and Dieter Sondermann.
\newblock Arrow's theorem, many agents, and invisible dictators.
\newblock \emph{Journal of Economic Theory}, 5\penalty0 (2):\penalty0 267--277,
  1972.

\bibitem[Lewis(1988)]{lewis1988infinite}
Alain~A Lewis.
\newblock An infinite version of arrow's theorem in the effective setting.
\newblock \emph{Mathematical Social Sciences}, 16\penalty0 (1):\penalty0
  41--48, 1988.

\bibitem[Mihara(1995)]{mihara1995arrow}
H~Reiju Mihara.
\newblock \emph{Arrow's Theorem, Turing computability, and oracles}.
\newblock PhD thesis, University of Minnesota, 1995.

\bibitem[Mihara(1997)]{mihara1997arrow}
H~Reiju Mihara.
\newblock Arrow's theorem and turing computability.
\newblock \emph{Economic Theory}, 10\penalty0 (2):\penalty0 257--276, 1997.

\bibitem[Mihara(1999)]{mihara1999arrow}
H~Reiju Mihara.
\newblock Arrow's theorem, countably many agents, and more visible invisible
  dictators.
\newblock \emph{Journal of Mathematical Economics}, 32\penalty0 (3):\penalty0
  267--287, 1999.

\bibitem[Morreau(2019)]{sep-arrows-theorem}
Michael Morreau.
\newblock {Arrow’s Theorem}.
\newblock In Edward~N. Zalta, editor, \emph{The {Stanford} Encyclopedia of
  Philosophy}. Metaphysics Research Lab, Stanford University, {W}inter 2019
  edition, 2019.

\bibitem[{Open Logic
  Project}(2022{\natexlab{a}})]{openlogicproject2022intermediate}
{Open Logic Project}.
\newblock {Intermediate Logic}, 2022{\natexlab{a}}.
\newblock URL
  \url{https://builds.openlogicproject.org/courses/intermediate-logic/il-screen.pdf}.

\bibitem[{Open Logic Project}(2022{\natexlab{b}})]{openlogicproject2022sample}
{Open Logic Project}.
\newblock {Sample Logic Text}, 2022{\natexlab{b}}.
\newblock URL
  \url{https://builds.openlogicproject.org/courses/sample/open-logic-sample.pdf}.

\bibitem[Soare(1999)]{soare1999recursively}
Robert~I Soare.
\newblock \emph{Recursively enumerable sets and degrees: A study of computable
  functions and computably generated sets}.
\newblock Springer Science \& Business Media, 1999.

\bibitem[Soare(2016)]{soare2016turing}
Robert~I Soare.
\newblock \emph{Turing computability: Theory and applications}.
\newblock Springer, 2016.

\end{thebibliography}

\end{document}